\documentclass[11pt,english]{smfart}

\usepackage[english]{babel}
\usepackage{amssymb,url,xspace,smfthm}
\usepackage{stmaryrd}
\usepackage{commandes}
\usepackage{xcolor}
\usepackage{dsfont}
\usepackage{graphicx}
\usepackage{tikz}
\usetikzlibrary{cd}
\usepackage{imakeidx}
\usepackage{enumerate}
\usepackage{algpseudocode}
\usepackage{algorithm}
\usepackage{subcaption}

\title[Horn inequalities on a quiver with an involution]{Horn inequalities on a quiver with an involution}
\alttitle{Inégalités de Horn sur un carquois avec une involution}
\date{8 January 2026}
\author[A. Médoc]{Antoine Médoc}
\address{IMAG, University of Montpellier, CNRS, Montpellier, France}
\email{antoine.medoc@umontpellier.fr} 
\urladdr{}
\thanks{I would like to thank my doctoral advisor Paul-\'Emile Paradan for his ideas and constant help during the writing of this paper.}
\keywords{semi-invariant, quiver involutions}

\makeindex[columns=2, title=Index, 
options= -s indexstyle.ist]

\begin{document}
\def\smfbyname{}

\begin{abstract}
	Derksen and Weyman described the cone of semi-invariants associated with a quiver \cite{derksen_semi-invariants_2000}.
	We give an inductive description of this cone, followed by an example of refinement of the inequalities characterising anti-invariant weights in the case of a quiver equipped with an involution.
\end{abstract}

\begin{altabstract}
	Derksen et Weyman ont décrit le cône des semi-invariants associés à un carquois \cite{derksen_semi-invariants_2000}.
	Nous donnons une description par récurrence de ce cône puis un exemple de raffinement des inégalités caractérisant les poids anti-invariants dans le cas d'un carquois muni d'une involution.
\end{altabstract}

\maketitle

\tableofcontents

\mainmatter

\section{Introduction}

A \textit{quiver} is a finite oriented graph and is described by a finite set of vertices $Q_0$, a finite set of arrows $Q_1$ and two maps maps $h,t$ from $Q_1$ to $Q_0$: any arrow $a\in Q_1$ has initial vertex $ta\in Q_0$ (its \textit{tail}) and terminal one $ha\in Q_0$ (its \textit{head}).
We will only consider quivers without oriented cycle.

For all $m,n\in\N^*$, we denote by $\Mat_{m,n}(\C)$ (resp. $\GL_n(\C)$) the set of all $m$ by $n$ matrices with complex coefficients (resp. of all invertible square matrices of order $n$ with complex coefficients).
A \textit{dimension} on $Q$ is a map $\alpha\in \N^{Q_0}$. 
The set of \textit{representations} of $Q$ of dimension $\alpha$ is
$$\Rep_\alpha(Q)= \bigoplus_{a\in Q_1} \Mat_{\alpha(ha),\alpha(ta)}(\C)$$ 
and a representation of $Q$ of dimension $\alpha$ is an element $V\in\Rep_\alpha(Q)$.
A subrepresentation of dimension $\beta\in \N^{Q_0}$ of $V\in \Rep_\alpha(Q)$ is a set $(W(x))_{x\in Q_0}$ such that, for all vertex $x$, $W(x)$ is a $\beta(x)$-dimensional linear subspace of $\C^{\alpha(x)}$ and, for all arrow $a$, $V(a) W(ta) \subset W(ha)$.\\  

\begin{defi}
	We denote $\beta\leqslant \alpha$ if, for all $x\in Q_0$, $\beta(x)\leqslant \alpha(x)$.
	We say that $\beta\in\N^{Q_0}$ is a \textit{generic subdimension} of $\alpha\in \N^{Q_0}$ if any representation of dimension $\alpha$ admits a subrepresentation of dimension $\beta$; in this case we denote $\beta\hookrightarrow \alpha$.
\end{defi}
Remark that $\beta\hookrightarrow \alpha\Rightarrow \beta\leqslant \alpha$. 
We denote by $\GL_\alpha(Q)$ the reductive algebraic group $\prod_{x\in Q_0} \GL_{\alpha(x)}(\C)$.
The natural action of $\GL_\alpha(Q)$ on $\Rep_\alpha(Q)$ yields an action on the ring of polynomial functions (or coordinate ring) $\C[\Rep_\alpha(Q)]$: for all $f\in \C[\Rep_\alpha(Q)]$, $g\in \GL_\alpha(Q)$ and $V\in\Rep_\alpha(Q)$,
\begin{align*}
	g\cdot V &= (g_{ha} V(a) g_{ta}^{-1} )_a,\\ 
	(g\cdot f)(V) &= f(g^{-1}\cdot V).
\end{align*}
A \textit{weight} $\s$ on $Q$ is an element of $\Z^{Q_0}$. To any weight $\s$ on $Q$ we associate a multiplicative character 
$$\appli{\chi_\s}{g}{\GL_\alpha(Q)}{\prod_{x\in Q_0} \det(g_x)^{\s(x)}}{\C^*}$$
and a ring morphism $\Z^{Q_0}\rightarrow \Z$ which associate $\alpha\in \Z^{Q_0}$ to the sum 
$$\s(\alpha) = \sum_{x\in Q_0} \s(x) \alpha(x).$$
A polynomial function $f\in \C[\Rep_\alpha(Q)]$ is a $\GL_{\alpha}(Q)$-\textit{semi-invariant} of weight $\s$ if, for all $g\in \GL_\alpha(Q)$, 
$g\cdot f = \chi_\s(g)f.$
We denote by $\SI(Q,\alpha)_\s$ the $\C$-linear subspace of all $\GL_{\alpha}(Q)$-semi-invariant of weight $\s$.
We are interested in the semigroup 
$$\Sigma(Q,\alpha)= \ens{ \s\in \Z^{Q_0} | \SI(Q,\ga)_\s \neq \ens 0 }.$$
In 1994, Alistair King showed that an element of this semi-group satisfies one homogeneous linear equation and a finite number of homogeneous linear inequalities \cite[Theorem 4.1]{king_moduli_1994}.
In 2000, Harm Derksen and Jerzy Weyman showed that a weight satisfying King's conditions is an element of the semi-group \cite[Theorem 3]{derksen_semi-invariants_2000}.\\

\begin{theo}[Derksen-Weyman] \label{theo-dw-sigma-is-a-cone} 
	For all weight $\s\in\Z^{Q_0}$, we have $\s\in \Sigma(Q,\alpha)$ if and only if the two following conditions are satisfied. 
	\begin{enumerate}
		\item $\s(\alpha)=0$;
		\item $\s(\beta)\leqslant 0$ for all $\beta \hookrightarrow\alpha$.
	\end{enumerate}
\end{theo}
The two authors then gave a full description of the faces of this cone \cite[Theorem 5.1]{derksen_combinatorics_2011}; Nicolas Ressayre later provided another proof of this description \cite[Theorem 3]{ressayre_git-cones_2012}.

By Theorem \ref{theo-dw-sigma-is-a-cone}, the semigroup $\Sigma(Q,\alpha)$ is saturated in the latice $\Z^{Q_0}$. Applied to the triple flag quiver \cite[Section 3]{derksen_semi-invariants_2000}, this general result on quivers yields the saturation of the Littlewood-Richardson coefficients obtained by Allen Knutson and Terence Tao \cite{knutson_honeycomb_1999}.  

The \textit{Euler-Ringel form} $\pscal\cdot\cdot$ associated to the quiver $Q$ is defined as follows: for all $\alpha,\beta\in\Z^{Q_0}$,
$$\pscal{\alpha}{\beta} = \sum_{x\in Q_0} \alpha(x)\beta(x) - \sum_{a\in Q_1} \alpha(ta) \beta(ha) \in\Z.$$
For all $\alpha\in \Z^{Q_0}$, we denote by $\pscal{\alpha}{\cdot}$ the weight associated to $\beta\mapsto \pscal{\alpha}{\beta}$.
The first result of this paper is an inductive description of the semigroup $\Sigma(Q,\alpha)$, in the spirit of Horn's conjecture \cite{horn_eigenvalues_1962,fulton_expose_1998}. 
It is a refinement of Theorem \ref{theo-dw-sigma-is-a-cone} and a consequence of Theorem \ref{theo-DW-no-involution-circ-1}.\\

\begin{theo} \label{theo-inductive-sigma}
	For all weight $\s\in\Z^{Q_0}$, $\s\in \Sigma(Q,\alpha)$ if and only if the two following conditions are satisfied. 
	\begin{enumerate}
		\item $\s(\alpha)=0$;
		\item $\s(\beta)\leqslant 0$ for all $\beta\leqslant \alpha$ such that $\pscal\beta\cdot \in \Sigma(Q,\alpha-\beta)$.
	\end{enumerate}
\end{theo}

In Section \ref{sec-inductive-sigma}, we provide a proof of Theorem \ref{theo-inductive-sigma} inspired by the techniques developed by Prakash Belkale in \cite{belkale_geometric_2006} (see also \cite{berline_horn_2018}).
This proof is interesting because it naturally adapts to the specific case of quivers equipped with an involution. In Section \ref{section-quiver-with-involutions}, we use the techniques introduced in Section \ref{sec-inductive-sigma} to prove Theorem \ref{theo-coro-main}, which concerns this particular case and constitutes the main result of this paper.\\

\begin{defi} \label{defi-involution}
	An \textit{involution} $\tau$ of a quiver $Q$ is a couple $(\tau_0,\tau_1)$ such that $\tau_0: Q_0\rightarrow Q_0$ and $\tau_1:Q_1\rightarrow Q_1$ are self-inverse maps satisfying $h\circ \tau_1 = \tau_0 \circ t$ and $t\circ \tau_1 = \tau_0 \circ h$. 
\end{defi}

See Example \ref{exemquiver1} below. We give other examples of quivers that admit an involution in Section \ref{sec-examples-of-quiver-with-involution}. See also \cite[Proposition 4.2]{derksen_generalized_2002}. For all matrix $M \in \Mat_{m,n}(\C)$, we denote by $M^*$ its conjugate transpose in $\Mat_{n,m}(\C)$.
Let $\tau$ be an involution on $Q$. 
It yields:
\begin{itemize}
	\item an involution $\tau$ on the ring $\Z^{Q_0}$ by $\tp \alpha = (\alpha(\tau x))_x$;
	\item an anti-linear isomorphism $\Rep_\alpha(Q) \rightarrow \Rep_{\tp \alpha}(Q)$ by 
	$$\tp V = (V(\tp a)^*)_a;$$
	\item a group isomorphism $\GL_\alpha(Q) \rightarrow \GL_{\tp \alpha}(Q)$ by $\tp g = (g(\tp x)^{*,-1})_x$;
	\item an anti-linear algebra isomorphism $\C[\Rep(Q,\alpha)] \simeq \C[\Rep(Q,\tau\cdot \alpha)]$
	by 
	$$(\tau\cdot f) (V) = \overline{f(\tau\cdot V)}.$$
	which restricts to $\SI(Q,\alpha)_\s \simeq\SI(Q,\tau\cdot \alpha)_{-\tau \cdot \s}$.
\end{itemize}
Then, if $\alpha=\tp \alpha$, $\tau$ yields an involution on $\Rep_\alpha(Q)$, $\GL_\alpha(Q)$ and 
$$\C[\Rep(Q,\alpha)] = \bigoplus_{\s\in \Z^{Q_0}} \SI(Q,\alpha)_\s$$  while $-\tau$ yields an involution on the saturated semigroup $\Sigma(Q,\alpha)$.\\

\begin{defi} \label{defi-symmetric-representation}
	A \textit{$\tau$-symmetric representation} of $Q$ is a representation $V$ of $Q$ such that $V=\tp V$, i.e. for all $a\in Q_1$, $V(a) = V(\tau a)^*$. 
\end{defi}

In particular, if $V$ is a $\tau$-symmetric representation then $\dim V = \tp \dim V$.
Our definition of a symmetric representation differs from the orthogonal and symplectic representations studied by Harm Derksen and Jerzy Weyman \cite[Section 2]{derksen_generalized_2002} because the involution on $\Rep_\alpha(Q)$ is anti-linear.

The main result of this paper is Theorem \ref{theo-coro-main} and it shows that, if a dimension satisfies $\alpha= \tp \alpha$, then the elements $\s\in \Sigma(Q,\alpha)$ satisfying $\s = -\tp \s$ are caracterised by a smaller set of inequalities than in Thereom \ref{theo-inductive-sigma}. It is a consequence of Theorem \ref{theo-symmetric-sigma-inequalities-circ1}.\\

\begin{theo} \label{theo-coro-main}
	Let $\alpha\in \N^{Q_0}$ and $\s\in \Z^{Q_0}$ such that $\alpha=\tp\alpha$ and $\s= -\tp \s$.
	Then $\s\in \Sigma(Q,\alpha)$ if and only if $\s(\beta)\leqslant 0$ for all $\beta\in\N^{Q_0}$ such that
	\begin{enumerate}
		\item $\ga:= \alpha-\beta-\tp\beta \in\N^{Q_0}$;
		\item $\pscal{\beta}{\cdot} \in \Sigma(Q,\ga) \cap \Sigma(Q,\tp\beta)$.
	\end{enumerate}
\end{theo}

The conditions of theorems \ref{theo-dw-sigma-is-a-cone}, \ref{theo-inductive-sigma} and \ref{theo-coro-main} are arithmetic conditions that allow us to give explicit examples of inequalities characterising the cone: see Section \ref{section-examples}. Example \ref{exemquiver1} below illustrates the reduction in the number of inequalities to characterise the anti-invariant elements of the semigroup. \\

\begin{exem} \label{exemquiver1}
	Let $Q_0$ be $\ens{x_i; 1\leqslant i \leqslant 6}$ a set of six vertices, $Q_1$ be $\ens{a_i; 1\leqslant i \leqslant 5}$ a set of five arrows and $t,h:Q_1\rightarrow Q_0$ two maps defined by $ta_1 = x_1$, $t a_2 = x_2$, $ha_1 = ha_2 = ta_3 = x_3$, $ha_3 = ta_4 = ta_5 = x_4$, $ha_4 = x_5$ and $ha_5 = x_6$. The permutation defined by $\tau_0 = (x_1 \; x_5) (x_3\; x_4) (x_2\; x_6)$ and $\tau_1:= (a_1 \; a_4) (a_2\; a_5)$ in cycle notation is an involution of the quiver $(Q_0,Q_1, h,t)$ (it is the only one up to quiver isomorphism).
	See figure \ref{figquiver2}.
	
	Let $\alpha$ a dimension vector such that $(\alpha(x_i))_{1\leqslant i \leqslant 6} = (2,3,4,4,3,2)$. We have $\alpha =\tp \alpha$.
	Theorem \ref{theo-dw-sigma-is-a-cone} characterises $\Sigma(Q,\alpha)$ with $244$ equations, Theorem \ref{theo-inductive-sigma} with $57$ equations.
	Map $\s \mapsto (\s(x_4), \s(x_5), \s(x_6))$ is a ring isomorphism from $\ens{ \s\in \Z^{Q_0} | \s= -\tp \s}$ to $\Z^{3}$. Theorem \ref{theo-coro-main} characterises $\Sigma(Q,\alpha)$ with the following $9$ equations (we omit $\s(0) \leqslant 0$): for all $\s\in \Z^{Q_0}$ such that $\s= -\tp\s$, $\s\in \Sigma(Q,\alpha)$ if and only if 
	\begin{align*}
		\s(x_6) &\leqslant 0,
		& \s(x_5) \leqslant 0 \\
		3\s(x_5) + 2\s(x_6)& \leqslant 0,
		&\s(x_4) + \s(x_6) \leqslant 0 \\
		\s(x_4) + 2\s(x_6) & \leqslant 0,
		&\s(x_4) + \s(x_5)  \leqslant 0 \\
		2\s(x_4) + 3\s(x_5) & \leqslant 0,
		&3\s(x_4) + 2\s(x_5) + \s(x_6) \leqslant 0 \\
		4\s(x_4) + 3\s(x_5) + 2\s(x_6) & \leqslant 0. &
	\end{align*}
	These inequalities are not independant: $\s(x_6) \leqslant 0$ and  $\s(x_5) \leqslant 0$ implies $3\s(x_5) + 2\s(x_6) \leqslant 0$.
	
	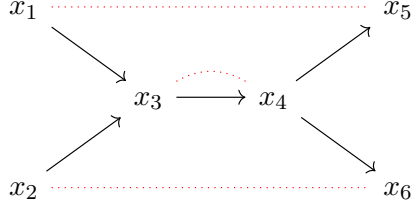
\begin{figure}[htbp]
		\centering
		\caption{A quiver equipped with an involution.}
		\label{figquiver2}
		\[\begin{tikzcd}
			{x_1} &&& {x_5} \\
			& {x_3} & {x_4} \\
			{x_2} &&& {x_6}
			\arrow[color={red}, dotted, no head, from=1-1, to=1-4]
			\arrow[from=1-1, to=2-2]
			\arrow[from=2-2, to=2-3]
			\arrow[color={red}, bend left, dotted, no head, from=2-2, to=2-3]
			\arrow[from=2-3, to=1-4]
			\arrow[from=2-3, to=3-4]
			\arrow[from=3-1, to=2-2]
			\arrow[color={red}, dotted, no head, from=3-1, to=3-4]
		\end{tikzcd}\]
	\end{figure}
\end{exem}

\section{An inductive description of the cone} \label{sec-inductive-sigma}

\subsection*{Notations}

We will mostly use notations from the book \cite{derksen_introduction_2017}. 
The set of all non-negative integers is denoted by $\N$, the set of all positive integers by $\N^*$. For all $i\in\N$ we denote $[i] = \ens{j \in \N | 1\leqslant j \leqslant i}$.
If the context is not ambiguous, a family $(\alpha_i)_{i\in I}$ indexed by a set $I$ will be referred to only as $(\alpha_i)$.
We fix a quiver $Q:= (Q_0,Q_1,h,t)$ without oriented cycles.
For all $\alpha\in\Z^{Q_0}$, $\pscal{\alpha}{\cdot}$ is the weight associated to $\beta\mapsto \pscal{\alpha}{\beta}$; $\pscal{\cdot}{\alpha}$ is the weight associated to $\beta\mapsto \pscal{\beta}{\alpha}$.
For all $k,n\in\N$, we denote by $\Gr(k,n)$ the projective variety of all $k$-dimensional linear subspaces of $\C^n$.
For all $\alpha,\beta \in \N^{Q_0}$ such that $\beta \leqslant \alpha$ we consider the projective variety
$\Gr(\beta, \alpha) = \prod_{x\in Q_0} \Gr(\beta(x),\alpha(x)).$

Thereafter, $\alpha,\beta,\ga,\delta\in\N^{Q_0}$ will denote dimensions on the quiver $Q$; $i,j,k,l,m,n\in\N$ will denote non-negative integers; $V,W$ will denote representations of the quiver $Q$; $\UU,\VV,\WW$ will denote sets of representations of $Q$; $\s\in\Z^{Q_0}$ will denote a weight on $Q$.
From section \ref{section-quiver-with-involutions} onwards, the quiver $Q$ admits an involution $\tau$.

\subsection{Homological dimensions}

Let $\alpha,\beta \in\N^{Q_0}$. To two representations $V\in\Rep_\alpha(Q)$ and $W\in \Rep_\beta(Q)$ we associate the linear map
$$\appli{d_W^V}{ \ph }{ \prod_{x\in Q_0} \Mat_{\beta(x),\alpha(x)}(\C) }{ (\ph(ha) V(a) - W(a) \ph(ta) }{ \prod_{a\in Q_1} \Mat_{\beta(ha),\alpha(ta)}(\C) };$$
the set $\Hom_Q(V,W)$ of all \textit{morphisms} $\ph:V\rightarrow W$ is defined as the kernel of $d_W^V$; the (Yoneda) \textit{extension group} of $V$ and $W$ is the cokernel of $d_V^W$ and denoted by $\Ext_Q(V,W)$ (it can be identified with the class of equivalences classes of extensions of $V$ by $W$).
Both are finite dimensional linear spaces and we denote their minimal value by
\begin{align*}
	\hom_Q(\alpha,\beta) &= \min_{(V,W)\in\Rep_\alpha(Q) \x \Rep_\beta(Q)} \dim \Hom_Q(V,W), \\
	\ext_Q(\alpha,\beta) &= \min_{(V,W)\in\Rep_\alpha(Q) \x \Rep_\beta(Q)} \dim \Ext_Q(V,W).
\end{align*}
Since the map $\Hom_Q(\cdot,\cdot)$ is upper semicontinuous for the Zariski topology on $\Rep_\alpha(Q) \x \Rep_\beta(Q)$, the integer $\hom_Q(\alpha,\beta)$ (resp. $\ext_Q(\alpha,\beta)$) is the generic value of $\dim \Hom_Q(\cdot,\cdot)$ (resp. $\dim \Ext_Q(\cdot,\cdot)$) on $\Rep_\alpha(Q)\x \Rep_\beta(Q)$.

An \textit{isomorphism} is a morphism $\ph$ such that, for all vertex $x$, $\ph(x)$ is bijective.
Remark that the $\GL_\alpha(Q)$-orbit of a representation $V$ is the set of all the representations isomorphic to $V$.
The \textit{kernel} of a morphism $\ph$ is $(\Ker \ph(x))_{x\in Q_0}$ and its \textit{image} is $(\Im \ph(x))_{x\in Q_0}$.

Remark that, if $\alpha,\beta\in\N^{Q_0}$, for all $V\in\Rep_\alpha(Q)$ and $W\in\Rep_\beta(Q)$, $\pscal \alpha\beta$ is the difference between the dimension of the codomain of $d_W^V$ and the dimension of its domain, i.e. 
$$\pscal{\alpha}{\beta} = \dim \Hom_Q(V,W) - \dim \Ext_Q(V,W)$$ 
by the rank theorem; from this we deduce 
\begin{equation} \label{eq-ringel=hom-ext}
	\pscal{\alpha}{\beta} = \hom_Q(\alpha,\beta) - \ext_Q(\alpha,\beta).
\end{equation}

\subsection{Schofield's polynomials}   \label{sec-DW-theorem-on-Schofield-polynomials}

Harm Derksen and Jerzy Weyman showed in \cite[Theorem 1]{derksen_semi-invariants_2000} that Schofield's polynomials \cite[Section 1]{schofield_semi-invariants_1991} span the finite-dimensional vector space of semi-invariants $\SI(Q,\ga)_\s$ for a given dimension $\ga$ and a given weight $\s$.
From this we deduce the Reciprocity Property \cite[Corollary 1]{derksen_semi-invariants_2000}: for all $\alpha,\beta\in\N^{Q_0}$,
\begin{equation} \label{eq-reciprocity-property}
	\dim \SI(Q,\beta)_{\pscal{\alpha}{\cdot}} = \dim \SI(Q,\alpha)_{-\pscal{\cdot}{\beta}}
\end{equation}
and this motivate the following definition. Let $\alpha,\beta\in\N^{Q_0}$.

\begin{defi} \label{defi-a-circ-b}
	We denote $\alpha\circ \beta = \dim (\SI(Q,\beta)_{\pscal{\alpha}{\cdot}}) \in\N$. 
\end{defi}

As stated instated in \cite[Remark 4]{derksen_combinatorics_2011}, another corollary of \cite[Theorem 1]{derksen_semi-invariants_2000} is the following.

\begin{coro} \label{coro-dw-circ}
	The following assertions are equivalent:
	\begin{enumerate}[$(i)$]
		\item $\alpha\circ \beta \neq 0$;
		\item $\hom_Q(\alpha,\beta) = \ext_Q(\alpha,\beta) = 0$;
		\item $\pscal{\alpha}{\cdot} \in\Sigma(Q,\beta)$;
		\item $-\pscal{\cdot}{\beta}\in \Sigma(Q,\alpha)$.
	\end{enumerate}
\end{coro}

\subsection{Generic filtration dimensions}

Let $\alpha\in\N^{Q_0}$ and $r\geqslant 2$.

\subsubsection{Definition}

For all representation $V\in\Rep_\alpha(Q)$, a \textit{filtration} of $V$ of length $r-1\in\N^*$ is a sequence $F = (F_i)_{i\in [r-1]}$ of subrepresentations of $V$ such that, for all $i\in [r-2]$, $F_i\subset F_{i+1}$; if we denote $F_{r}=V$ and $F_0=\ens0$, the \textit{dimension} of the filtration $F$ is the sequence $(\dim F_i-\dim F_{i-1})_{i} \in (\N^{Q_0})^{r}$.
Let $(\alpha_i)\in (\N^{Q_0})^r$ be a sequence of $r$ dimensions such that 
$$\alpha = \sum_{i=1}^{r} \alpha_i.$$ 
We denote by $\Rep_{(\alpha_i)}(Q)$ the set of all $V \in \Rep_\alpha(Q)$ such that there exists a filtration of $V$ of dimension $(\alpha_i)$. 
For example, if $Q_1=\emptyset$ then $\Rep_{(\alpha_i)}(Q) = \Rep_\alpha(Q)$.
The following lemma is \cite[Lemma 3.5]{derksen_combinatorics_2011}.\\

\begin{lemm} \label{lemm-filtration-zariski-closed} 
	The set $\Rep_{(\alpha_i)}(Q)$ is Zariski-closed in $\Rep_\alpha(Q)$.
\end{lemm}

\begin{proof} 
	We adapt the proof of \cite[Lemma 3.1]{schofield_general_1992}, which is the case $r=2$.
	For all $i\in [r]$ we denote $\beta_i = \sum_{j\in [i]} \alpha_j$.
	We consider the set $Z$ of all $(V,(R_i))\in \Rep_\alpha(Q) \x \prod_{i\in [r-1]} \Gr(\beta_i,\beta_{i+1})$ such that $R$ is a filtration of $V$.
	It is a subvariety of $\Rep_\alpha(Q) \x \prod_{i\in [r-1]} \Gr(\beta_i,\beta_{i+1})$ hence its image $\Rep_{(\alpha_i)}(Q)$ by the projection on $\Rep_\alpha(Q)$ is closed.
\end{proof}

\begin{defi}
	The sequence $(\alpha_i)$ is a \textit{generic filtration dimension} if $\Rep_{(\alpha_i)}(Q) = \Rep_\alpha(Q).$
\end{defi}

\subsubsection{Schofield's theorem}

Remark that, for all $\beta\leqslant \alpha$, $(\beta,\alpha-\beta)$ is a generic filtration dimension if and only if $\beta\hookrightarrow \alpha$. 
For example, Lemma \ref{lemm-filtration-zariski-closed} shows that any subrepresentation $W$ of a general representation $V\in \Rep_\alpha(Q)$ satisfies $\dim W \hookrightarrow \alpha$.
It is a transitive relation: for all $\beta,\ga\in\N^{Q_0}$, if $\beta\hookrightarrow \alpha$ and $\ga \hookrightarrow \beta$ then $\ga\hookrightarrow \alpha$.
The following characterisation is Schofield's Theorem \cite[Theorem 3.3]{schofield_general_1992}: for all $\beta\in \N^{Q_0}$,
\begin{equation} \label{theo-schofield-subdimension-carac}
	\alpha \hookrightarrow \alpha+\beta \Leftrightarrow \ext_Q(\alpha,\beta)=0
\end{equation}

\subsubsection{Ressayre's dominant morphisms}

Being a generic filtration dimension corresponds to the dominant feature of a certain morphism: this will yield a useful inequality on the dimension of its domain and codomain.
We are working in the same context as in \cite[Section 4]{ressayre_multiplicative_2011}, which is a particular case of \cite[Section 2.1]{ressayre_multiplicative_2011}.
Let $(\alpha_i)\in (\N^{Q_0})^r$ be a sequence of $r$ dimensions such that $$\alpha = \sum_{i\in[r]} \alpha_i.$$
We denote $\beta_0= (0)_{x} \in \N^{Q_0}$ and, for all $i\in [r]$, $\beta_i = \sum_{j\in[i]} \alpha_j$.
We define a one-parameter subgroup $\la:\C^*\rightarrow \GL_\alpha(Q)$ by: for all $x\in Q_0$ and $t\in\C^*$,
$$ (\la)_x(t) = \diag( t^i \Id_{\alpha_i(x)} ; i\in[r] ) \in \GL_{\alpha(x)}(\C).$$
For all $x\in Q_0$, we dentote by $(e_j(x))_j$ the canonical basis of $\C^{\alpha(x)}$. For all $i\in [r]$ we denote
\begin{align*}
	E_i^\lambda&= \ens{ v\in \prod_{x\in Q_0} \C^{\alpha(x)} | \forall t\in \C^*, \lambda(t) v = t^i v }
\end{align*}
i.e. for all $x\in Q_0 $, $E_i^\lambda(x) = \spa \ens{ e_j(x) ; 1+\beta_{i-1}(x) \leqslant j \leqslant \beta_i(x) }$. We have
$$ \prod_{x\in Q_0} \C^{\alpha(x)} = \bigoplus_{i\in [r]} E_i^\lambda.$$
and, for all $i\in[r]$, $\alpha_i = \dim E_i^\lambda$.
Consider the parabolic subgroup of $\GL_\alpha(Q)$ associated to $\lambda$
$$ P(\lambda):= \ens{ p\in \GL_\alpha(Q) | \exists p_0\in\GL_\alpha(Q), \lambda(t) p \lambda(t)^{-1} \underset{t\rightarrow +\infty}{\longrightarrow} p_0 },$$
the set of fixed representation
$$ C^0(\lambda):= \ens{ V\in \Rep_\alpha(Q) | \forall t\in \C^*, \lambda(t)\cdot V = V } $$
and the associated Bia\l{}ynicki-Birula cell \cite{bialynicki-birula_theorems_1973}
$$C^+(\lambda) = \ens{ V\in \Rep_\alpha(Q) | \exists V_0\in C^0(\la), \lambda(t)\cdot V \underset{t\rightarrow +\infty}{\longrightarrow} V_0 }.$$
These last two sets, $C^0(\lambda)$ and $C^+(\lambda)$, are two linear subspaces of $\Rep_\alpha(Q)$ and satisfy $C^0(\lambda) \subset C^+(\lambda)$.
The cell $C^+(\lambda)$ is stable by the action of the subgroup $P(\lambda)$.
From this we know that the group $\GL_\alpha(Q)\x P(\lambda)$ acts on $\GL_\alpha(Q) \x C^+(\lambda)$ in the following way: for all $(g,p)\in \GL_\alpha(Q)\x P(\lambda)$ and $(g',V)\in \GL_\alpha(Q) \x C^+(\lambda)$,
$$(p,g)\cdot(g',V) = (gg'p^{-1}, p\cdot V).$$
We denote $\GL_\alpha(Q) \x_{P(\lambda)} C^+(\lambda)$ the quotient of $\GL_\alpha(Q) \x C^+(\lambda)$ by the action of $P(\lambda)$. The natural $P(\lambda)$-invariant morphism
$$\appli{\rho_\lambda}{(g,V)}{\GL_\alpha(Q)\x C^+(\lambda)}{ g\cdot V }{\Rep_\alpha(Q)}$$
yields a $\GL_\alpha(Q)$-equivariant morphism $\eta_\lambda:\GL_\alpha(Q)\x_{P(\lambda)} C^+(\lambda) \rightarrow \Rep_\alpha(Q)$.
For all $j\in[r]$ we denote 
$$F_j^\lambda:= \bigoplus_{i\in[j]} E_i^\lambda.$$
We have
\begin{align}
	C^0(\lambda) &= \ens{ V\in \Rep_\alpha(Q) | \forall i\in[r], E_i^\lambda <V}, \notag\\
	C^+(\lambda) &= \ens{V \in \Rep_\alpha(Q) | \forall j \in[r], F_j^\lambda <V}, \label{eq-c+lambda} \\
	P(\lambda) &= \ens{ g\in \GL_\alpha(Q) | \forall j \in[r], g(F_j^\lambda) < F_j^\lambda }, \label{eq-plambda}\\
	\Im (\eta_\lambda) &= \Im(\rho_\lambda) = \Rep_{(\alpha_i)}(Q). \label{eq-im-of-etalambda}
\end{align}
As in the proof of \cite[Theorem 3]{ressayre_multiplicative_2011} we compute $\delta(\lambda)$ the difference between the dimension of the codomain and of the domain of $\eta_\lambda$:
\begin{align*}
	\delta(\lambda) &= \dim \Rep_\alpha(Q) - \dim \GL_\alpha(Q) \x_{P(\lambda)} C^+(\lambda) \\
	&= \dim \Rep_\alpha(Q) -\dim \GL_\alpha(Q) + \dim P(\lambda) - \dim C^+(\lambda).
\end{align*}
with
\begin{align*}
	\dim \Rep_\alpha(Q) &= \sum_{a} \sum_{i,j} \alpha_i(ta) \alpha_j(ha), & \dim \GL_\alpha(Q) &= \sum_{x} \sum_{i,j} \alpha_i(x)\alpha_j(x)\\
	\intertext{and, by Equations \eqref{eq-c+lambda} and \eqref{eq-plambda},}
	\dim C^+(\lambda) &= \sum_{a} \sum_{j\leqslant i} \alpha_i(ta) \alpha_j(ha),&
	\dim P(\lambda) &= \sum_{x} \sum_{j\leqslant i} \alpha_i(x) \alpha_j(x)
\end{align*}
From this we deduce
$$\delta(\lambda) = -\sum_{i<j} \pscal{\alpha_i}{\alpha_j}.$$
Yet, by Equation \eqref{eq-im-of-etalambda}, $(\alpha_i)$ is a generic filtration dimension if and only if $\eta_\lambda$ is surjective (if and only if $\eta_\lambda$ is dominant by Lemma \ref{lemm-filtration-zariski-closed}). This yields the following lemma.\\

\begin{lemm} \label{lemm-sum-ringel-generic-dim-vector}
	If $(\alpha_i)$ is a generic filtration dimension then
	$$\sum_{i<j} \pscal{\alpha_i}{\alpha_j} \geqslant 0.$$
\end{lemm}

It can be compared with the codimension condition \cite[Equation 0.1]{belkale_geometric_2006} and will be used together with Lemma \ref{lemm-hom-minimum-witness-no-symmetry}.

\subsection{Number of subrepresentations}

Let $\alpha,\beta \in\N^{Q_0}$.
In \cite[Section 4.1]{ressayre_multiplicative_2011}, Nicolas Ressayre consider the following definition: if a general representation of dimension $\alpha+\beta$ admits a finite number of subrepresententations of dimension $\alpha$, we denote this number by $\alpha\circ' \beta$; otherwise, we denote $\alpha\circ' \beta =0$. 
From this point of view, $\alpha\circ' \beta$ is the cardinal of a general fiber.
The following theorem shows that if $\pscal{\alpha}{\beta}=0$ then $\alpha \circ \beta = \alpha\circ' \beta$.
Recall that if $\alpha\circ\beta\neq 0$ then $\pscal{\alpha}{\beta}=0$.
Theorem \ref{theo-dw-schofield-number-of-subrepresentations} is \cite[Theorem 1]{derksen_number_2007} and will be used together with Lemma \ref{lemm-unique-witness}.\\

\begin{theo}[Derksen-Schofield-Weyman] \label{theo-dw-schofield-number-of-subrepresentations} 
	If $\pscal\alpha\beta =0$ then a general representation $V\in \Rep_{\alpha+\beta}(Q)$ admits exactly $\alpha\circ \beta$ subrepresentations of dimension $\alpha$.
\end{theo}

In a certain way, a subrepresentation of a general representation is also in general position \cite[Proposition 13]{derksen_number_2007}.\\

\begin{prop} \label{prop-subrepresentations-in-general-position}
	Assume $\beta\leqslant \alpha$ and $\pscal{\beta}{\alpha-\beta} =0$. Let $\UU_\beta$ a $\GL_\beta$-stable non-empty Zariski-open subset of $\Rep_\beta(Q)$. There exists $\UU_\alpha$ a non-empty Zariski-open subset of $\Rep_\alpha(Q)$ such that, for all $V\in\UU_\alpha$, any subrepresentation of $V$ lies in $\UU_\beta$.
\end{prop}

From this and the last theorem follows Corollary \ref{coro-general-sub-sub-representations} about the generic number of subrepresentations of subrepresentations.\\

\begin{coro} \label{coro-general-sub-sub-representations}
	Let $\alpha_1,\alpha_2,\alpha_3\in \N^{Q_0}$ such that $$\pscal{\alpha_1}{\alpha_2} = \pscal{\alpha_1+\alpha_2}{\alpha_3} = 0.$$ There exists a non-empty open subset $\UU$ of $\Rep_{\alpha_1+\alpha_2+\alpha_3}(Q)$ such that, for all $V\in \UU$ and all subrepresentation $W<V$ of dimension $\alpha_1+\alpha_2$, $W$ admits exactly $\alpha_1\circ \alpha_2$ subrepresentations of dimension $\alpha_1$.
\end{coro}

\begin{proof}
	By Derksen-Schofield-Weyman's Theorem \ref{theo-dw-schofield-number-of-subrepresentations}, there exists a non-empty open subset $\UU'$ of $\Rep_{\alpha_1+\alpha_2}(Q)$ such that, for all $W\in\UU'$, $W$ admits exactly $\alpha_1\circ \alpha_2$ subrepresentations.
	The union of orbits $\GL_{\alpha_1+\alpha_2}(Q)\cdot \UU'$ is a non-empty open subset of $\Rep_{\alpha_1+\alpha_2}(Q)$ and is stable by $\GL_{\alpha_1+\alpha_2}(Q)$. We conclude by Proposition \ref{prop-subrepresentations-in-general-position}.
\end{proof}

We conclude this section with the following multiplicative formula obtained by Nicolas Ressayre \cite[Theorem B]{ressayre_multiplicative_2011} (see also \cite[Section 7.4]{derksen_combinatorics_2011}).\\

\begin{prop} \label{prop-ressayre-multiplicative-coro2} 
	Assume that $\pscal{\beta_1}{\beta_2} = \pscal{\beta_1}{\beta_3} = 0$ and $\beta_2 \circ \beta_3 = 1$. We have
	$\beta_1\circ (\beta_2+\beta_3) = (\beta_1\circ \beta_2) \x (\beta_1\circ \beta_3).$
\end{prop}

\subsection{Discrepancy}

The definition and notation $\disc(V,\s)$ below comes from \cite[Definition 29]{chindris_simultaneous_2021}.\\
\begin{defi} \label{defi-discrepancy}
	The \textit{discrepancy} of a representation $V$ with respect to a weight $\s\in \Z^{Q_0}$ is 
	$$\disc(V,\s) = \max_{W<V} \s(\dim W)$$
	and a \textit{witness} for $(V,\s)$ is a subrepresentation $W<V$ reaching this maximum: $\s(\dim W)=\disc(V,\s)$.
	The discrepancy of a dimension $\alpha\in\N^{Q_0}$ with respect to a weight $\s\in \Z^{Q_0}$ is
	$$\disc(\alpha,\s) = \max_{\beta\hookrightarrow\alpha} \s(\beta).$$
\end{defi}

Remark that, by Lemma \ref{lemm-filtration-zariski-closed}, there exists a non-empty Zariski-open subset $\VV\subset \Rep_\alpha(Q)$ such that, for all $V\in \VV$, $\disc(V,\s) = \disc(\alpha,\s)$. 
The next lemma is a Harder-Narasimhan type lemma and is stated in \cite[Proposition 3.1]{huszar_non-commutative_2021}. It can be compared with \cite[Lemma 5]{belkale_local_2001} \\

\begin{lemm} \label{lemm-unique-witness} 
	If $W_1,W_2$ are witnesses for $(V,\s)$, then so are $W_1\cap W_2$ and $W_1+W_2$. In particular, there exists a unique minimum (resp. maximum) witness with respect to the inclusion.
\end{lemm}

\begin{proof}
	Let $W_1,W_2$ be two witnesses for $(V,\s)$. Since $\disc(V,\s)$ is a maximum we have
	\begin{align*}
		\s(\dim (W_1+W_2)) &= \s(\dim W_1) + \s(\dim W_2) - \s(\dim (W_1\cap W_2)) \\
		&\geqslant 2\disc(V,\s) - \disc(V,\s) = \disc(V,\s),\\
		\intertext{and}
		\s(\dim (W_1\cap W_2)) &= \s(\dim W_1) + \s(\dim W_2) - \s(\dim (W_1+W_2)) \\
		& \geqslant 2\disc(V,\s) - \disc(V,\s) = \disc(V,\s),
	\end{align*}
	hence $\s (\dim (W_1+W_2)) = \s(\dim (W_1\cap W_2)) = \disc(V,\s)$. The first point is proven.
	We consider the dimension
	$$\pappli{d}{x}{Q_0}{\min_{W<V, \s(\dim W) = \disc(V,\s)} W(x) }{\N}.$$
	For all $x\in Q_0$ there exists $W_x<V$ such that $d(x) = \dim W_x(x)$. By the first point of the lemma, $m:= \bigcap_{x\in Q_0} W_x$ is a witness for $(V,\s)$ such that $\dim m \leqslant d$.  
	Let $W$ be a witness. For all $x\in Q_0$, since $m\cap W$ is a witness then $d(x) \leqslant \dim (W\cap m)(x)$. From this we know that $m\leqslant d \leqslant \dim m\cap W$ and $m = m\cap W$, i.e. $m< W$. 
	Similarly, there exists a maximum witness.
\end{proof}

\begin{lemm} \label{lemm-hom-minimum-witness-no-symmetry}
	Let $W$ be the minimum witness of $(V,\s)$ and $\beta= \dim W$ its dimension. We have $\hom_Q(\beta,\alpha-\beta)=0$.
\end{lemm}

\begin{proof}
	Let $\ph\in\Hom_Q(W,V/W)$, $i=\dim\Im\ph$ and $k:= \dim \Ker\ph$. Let $\pi:V\twoheadrightarrow V/W$ the canonical projection on the quotient $V/W$ and $\tilde I$ the subrepresentation $\pi^{-1}(\Im \ph) <V$. We have
	$$\disc(V,\s) \geqslant \s(\dim\tilde I) = \s(i)+\s(\beta) = \s(i)+\disc(V,\s)$$
	hence $\s(i)\leqslant 0$. Yet, by the rank theorem, $\beta = i+k$ and 
	$$\disc(V,\s) = \s(\beta) = \s(i)+\s(k) \leqslant \s(k)$$
	hence the kernel $\Ker\ph$ is a witness for $(V,\s)$ and a subrepresentation of the minimum witness $W$, i.e. $K=W$, i.e. $\ph=0$. From this we deduce $\hom_Q(\beta, \ga-\beta) =0$. 
\end{proof}

We conclude this section by showing that the maximum that defines the discrepancy of $\alpha$ with respect to $\s$ is reached on a limited set of dimensions.\\

\begin{lemm} \label{lemm-disc-calculation-circ1}
	We have
	$$\disc(\alpha,\s) = \max_{  \beta\circ (\alpha-\beta) =1 } \s(\beta).$$
\end{lemm}

\begin{proof}
	It suffices to show that there exists $\beta\leqslant \alpha$ such that $\beta\circ(\alpha-\beta) = 1$ and $\s(\beta)=\disc(\alpha,\s)$.
	By Lemma \ref{lemm-filtration-zariski-closed}, there exists a non-empty open subset $\UU$ of $\Rep_\alpha(Q)$ such that, for all $\beta$ which do not satisfy $\beta \hookrightarrow \alpha$ and for all $V\in \UU$, $V$ do not admit any subrepresentation of dimension $\beta$. 
	For all $\beta \leqslant \alpha$ such that $\pscal{\beta}{\alpha-\beta} =0$, by Proposition \ref{prop-subrepresentations-in-general-position}, there exists $\VV_\beta$ a non-empty open subset of $\Rep_\alpha(Q)$ such that, for all $V\in\VV_\beta$, $V$ has exactly $\beta \circ (\alpha-\beta)$ subepresentations of dimension $\beta$. 
	We condider the intersection 
	$$\VV:= \bigcap_{\pscal{\beta}{\alpha-\beta} =0} \VV_\beta.$$
	There exists $V\in \UU\cap \VV$. Let $W$ be the minimum witness of $(V,\s)$ and $\beta:= \dim W$ its dimension. Since $V\in\UU$, we know that $\disc(\alpha,\s) = \s(\beta)$ and that $\beta \hookrightarrow \alpha$. 
	By Schofield's Theorem \ref{theo-schofield-subdimension-carac}, $\ext_Q(\beta, \alpha-\beta)=0$.
	Yet, using lemma \ref{lemm-hom-minimum-witness-no-symmetry}, $\hom_Q(\beta,\alpha-\beta) = 0$. In particular, $\pscal{\beta}{\alpha-\beta} =0$ (we could also use Lemma \ref{lemm-sum-ringel-generic-dim-vector} instead of Theorem \ref{theo-schofield-subdimension-carac}).
	Yet, by Lemma \ref{lemm-unique-witness} we also know that $V$ admits a unique subrepresentation of dimension $\beta$ i.e. (since $V\in\VV$) $\beta \circ (\alpha-\beta) = 1$.
\end{proof}

\subsection{Schofield's Equation on discrepancy}

Let $\alpha,\beta\in \N^{Q_0}$ be two dimensions and $\s\in\Z^{Q_0}$ be a weight on $Q$. 
The following equation is \cite[Theorem 5.4]{schofield_general_1992}:
\begin{equation} \label{eq-schofield} 
	\ext_Q(\alpha,\beta) = \disc(\alpha, -\pscal{\cdot}{\beta}) = 
	\disc( \beta, \pscal{\alpha}{\cdot} ) - \pscal{\alpha}{\beta}.
\end{equation}
On one hand, together with Schofield's Theorem \ref{theo-schofield-subdimension-carac}, this gives an algorithm \cite[Section 5]{schofield_general_1992} to compute the integer $\ext_Q(\alpha,\beta)$ (hence $\hom_Q(\alpha,\beta)$) and the assertion $\beta\hookrightarrow \alpha$.
On the other hand, together with Corollary \ref{coro-dw-circ}, we deduce the Generalized Saturation Theorem \cite[Corollary 2.19]{derksen_combinatorics_2011}. In particular, the semigroup $\Sigma(Q,\alpha)$ is saturated in the lattice $\Z^{Q_0}$ \cite[Theorem 3]{derksen_semi-invariants_2000}.

\subsection{King's criterion on semi-stability}

Let $\alpha \in\N^{Q_0}$ and $\s\in \Z^{Q_0}$. 
The following definition corresponds to the $\GL_\alpha(Q)$-semi-stability of an element $(V,1) \in \Rep_\alpha(Q) \oplus \chi_\s$ in Geometric Invariant Theory.\\

\begin{defi}
	A representation $V$ is $\s$-\textit{semi-stable} if there exists $n\in\N$ and $f\in\SI(Q,\dim V)_{n\s}$ such that $f(V)\neq 0$.
\end{defi}

In particular, if $\s\in\Sigma(Q,\alpha)$ then there exists a $\s$-stable representation of dimension $\alpha$. By \cite[Theorem 1]{derksen_semi-invariants_2000}, the reciprocal is true. This will be used in Section \ref{sec-descriptions-of-the-cone} together with Lemma \ref{lemm-disc-calculation-circ1}.\\

\begin{lemm}[Derksen-Weyman] \label{lemm-sigma-made-of-semi-stable}
	We have $\s\in \Sigma(Q,\alpha)$ if and only if there exists a $\s$-semi-stable representation $V$ of dimension $\alpha$.
\end{lemm}

King's criterion \ref{theo-king-criterion} below is \cite[Theorem 4.1]{king_moduli_1994} and establishes a link between the $\s$-semi-stability of a representation $V$ and the discrepancy of the couple $(V,\s)$.\\

\begin{theo}[King's criterion]  \label{theo-king-criterion}
	The representation $V$ is $\s$-semi-stable if and only if $\s(\dim V) = \disc(V,\s) = 0$.
\end{theo}

\subsection{Horn type inequalities on the semi-invariants} \label{sec-descriptions-of-the-cone}
Let $\alpha\in\N^{Q_0}$.
Lemmas \ref{lemm-filtration-zariski-closed}, \ref{lemm-sigma-made-of-semi-stable} and King's criterion \ref{theo-king-criterion} yield Derksen-Weyman's Theorem \ref{theo-dw-sigma-is-a-cone}: for all $\s\in\Z^{Q_0}$, 
$$\s\in \Sigma(Q,\alpha) \Leftrightarrow \disc(\alpha,\s) = \s(\alpha)=0.$$
Derksen-Weyman's Theorem \ref{theo-dw-sigma-is-a-cone} and Lemma \ref{lemm-disc-calculation-circ1} yield Theorem \ref{theo-DW-no-involution-circ-1}; Theorem \ref{theo-DW-no-involution-circ-1} and Corollary \ref{coro-dw-circ} yield the inductive description of $\Sigma(Q,\alpha)$ stated in Theorem \ref{theo-inductive-sigma}. In fact, the proof of Lemma \ref{lemm-disc-calculation-circ1} shows that we do not need Theorem \ref{theo-dw-schofield-number-of-subrepresentations} to prove Theorem \ref{theo-inductive-sigma}.
\\

\begin{theo} \label{theo-DW-no-involution-circ-1}
	For all weight $\s\in\Z^{Q_0}$, $\s\in \Sigma(Q,\alpha)$ if and only if the two following conditions are satisfied:
	\begin{enumerate}
		\item $\s(\alpha)=0$;
		\item $\s(\beta) \leqslant 0$ for all $\beta\leqslant \alpha$ such that $\beta\circ (\alpha-\beta)=1$.
	\end{enumerate}
\end{theo}

Section \ref{sec-inductive-sigma} and the proof of this last theorem allowed us to introduce the useful tools for proving Theorem \ref{theo-coro-main} in the following section.

\section{Quivers with an involution} \label{section-quiver-with-involutions}

\subsection{Examples} \label{sec-examples-of-quiver-with-involution}

\begin{figure}[htbp]
	\centering
	\caption{The oriented straight quiver $A_n^\rightarrow$}
	\label{fig-quiver-an}
	\[\begin{tikzcd}
		1 & 2 & \cdots & {n-1} & n
		\arrow[from=1-1, to=1-2]
		\arrow[color=red, bend left=15, dotted, no head, from=1-1, to=1-5]
		\arrow[from=1-2, to=1-3]
		\arrow[color=red, bend right=15, dotted, no head, from=1-2, to=1-4]
		\arrow[from=1-3, to=1-4]
		\arrow[from=1-4, to=1-5]
	\end{tikzcd}\]
\end{figure}
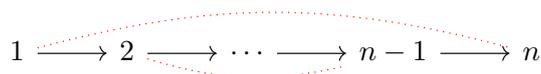

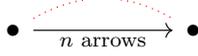
\begin{figure}[htbp]
	\centering
	\caption{The $n$-Kronecker quiver $\Theta_n$}
	\label{fig-quiver-kronecker}
	\[\begin{tikzcd}
		\bullet && \bullet
		\arrow["{n \text{ arrows}}"', from=1-1, to=1-3]
		\arrow[color=red, bend left, dotted, no head, from=1-1, to=1-3]
	\end{tikzcd}\]
\end{figure}

See also \cite[Proposition 4.2]{derksen_generalized_2002} and Example \ref{exemquiver1} for examples of quivers with an involution.
A quiver admits an involution if and only if its opposite quiver admits an involution. 
The oriented straight quiver $A_n^\rightarrow$ (see Figure \ref{fig-quiver-an}) and the $n$-Kronecker quiver $\Theta_n$ (see Figure \ref{fig-quiver-kronecker}) admit an involution. 

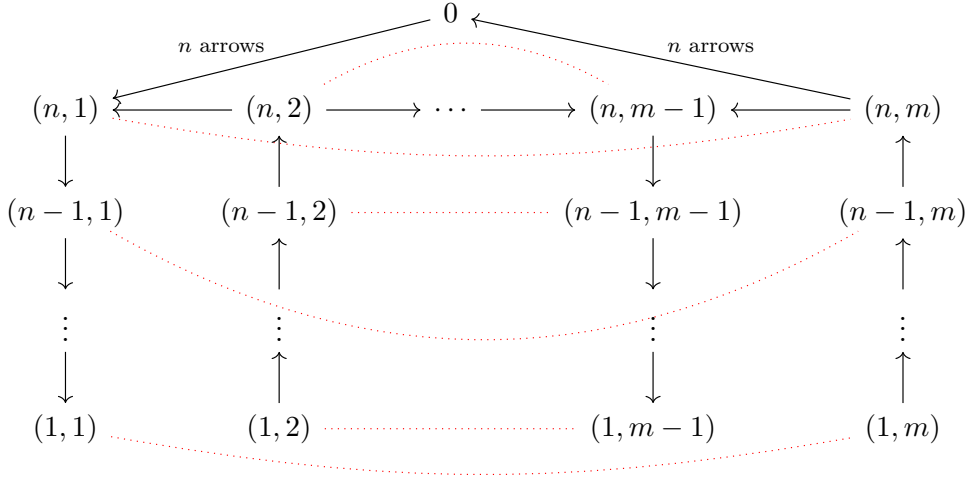
\begin{figure}[htbp]
	\centering
	\caption{The $(n,m)$-Chindris quiver if $m$ is even}
	\label{fig-quiver-chindris}
	\[\begin{tikzcd}
		&& 0 \\
		{(n,1)} & {(n,2)} & \cdots & {(n,m-1)} & {(n,m)} \\
		{(n-1,1)} & {(n-1,2)} && {(n-1,m-1)} & {(n-1,m)} \\
		\vdots & \vdots && \vdots & \vdots \\
		{(1,1)} & {(1,2)} && {(1,m-1)} & {(1,m)}
		\arrow["{n \text{ arrows}}"', from=1-3, to=2-1]
		\arrow[color=red, bend right=10, no head, dotted, from=2-1, to=2-5]
		\arrow[from=2-1, to=3-1]
		\arrow[from=2-2, to=2-1]
		\arrow[from=2-2, to=2-3]
		\arrow[color=red, bend left, dotted, no head, from=2-2, to=2-4]
		\arrow[from=2-3, to=2-4]
		\arrow[from=2-4, to=3-4]
		\arrow["{n \text{ arrows}}"', from=2-5, to=1-3]
		\arrow[from=2-5, to=2-4]
		\arrow[shift left, color=red, bend right, dotted, no head, from=3-1, to=3-5]
		\arrow[from=3-1, to=4-1]
		\arrow[from=3-2, to=2-2]
		\arrow[color=red, dotted, no head, from=3-2, to=3-4]
		\arrow[from=3-4, to=4-4]
		\arrow[from=3-5, to=2-5]
		\arrow[from=4-1, to=5-1]
		\arrow[from=4-2, to=3-2]
		\arrow[from=4-4, to=5-4]
		\arrow[from=4-5, to=3-5]
		\arrow[color=red, bend right=10, dotted, no head, from=5-1, to=5-5]
		\arrow[from=5-2, to=4-2]
		\arrow[color=red, dotted, no head, from=5-2, to=5-4]
		\arrow[from=5-5, to=4-5]
	\end{tikzcd}\]
\end{figure}

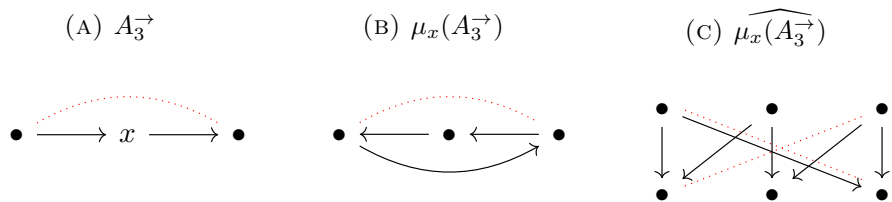
\begin{figure}[htbp]
	\centering
	\caption{Mutation and doubling of the oriented straight quiver $A_3^\rightarrow$}
	\label{fig-quiver-a3-mutated-doubling}
	\begin{subfigure}{0.23\textwidth}
		\caption{$A_3^\rightarrow$} 
		\label{subfig-a3}
		\[\begin{tikzcd}
			\bullet & x & \bullet
			\arrow[from=1-1, to=1-2]
			\arrow[color=red, bend left, dotted, no head, from=1-1, to=1-3]
			\arrow[from=1-2, to=1-3]
		\end{tikzcd}\]
	\end{subfigure}
	\hspace{1cm}
	\begin{subfigure}{0.23\textwidth}
		\caption{$\mu_x(A_3^\rightarrow)$} \label{subfig-mutated-a3}
		\[\begin{tikzcd}
			\bullet & \bullet & \bullet
			\arrow[color=red, bend left, dotted, no head, from=1-1, to=1-3]
			\arrow[bend right, from=1-1, to=1-3]
			\arrow[from=1-2, to=1-1]
			\arrow[from=1-3, to=1-2]
		\end{tikzcd}\]
	\end{subfigure}
	\hspace{1cm}
	\begin{subfigure}{0.23\textwidth}
		\caption{$\widehat{\mu_x(A_3^\rightarrow)}$} \label{subfig-doubling-mutated-a3}
		\[\begin{tikzcd}
			\bullet & \bullet & \bullet \\
			\bullet & \bullet & \bullet
			\arrow[from=1-1, to=2-1]
			\arrow[from=1-1, to=2-3]
			\arrow[shift left, color=red, dotted, no head, from=1-1, to=2-3]
			\arrow[from=1-2, to=2-1]
			\arrow[from=1-2, to=2-2]
			\arrow[color=red, dotted, no head, from=1-3, to=2-1]
			\arrow[from=1-3, to=2-2]
			\arrow[from=1-3, to=2-3]
		\end{tikzcd}\]
	\end{subfigure}
\end{figure}

The Sun quiver (see section \ref{subsection-sun-quiver}) admit an involution.
If $m$ is even, the quiver defined in \cite[Section 3]{chindris_quivers_2009} (see Figure \ref{fig-quiver-chindris}) admits an involution.
We briefly mention the case of a quiver with oriented cycle in Figure \ref{fig-quiver-a3-mutated-doubling}. If a quiver $Q$ without oriented cycle admits an involution (see Figure \ref{subfig-a3}), its mutated quiver $\mu_x(Q)$ \cite[page 9]{fomin_total_2010} (see Figure \ref{subfig-mutated-a3}) with respect to a fixed vertex $x$ also admits an involution (yet it may have an oriented cycle). If $Q$ (possibly with an oriented cycle) admits an involution, $\hat Q$ defined in \cite[Section 6.2]{derksen_combinatorics_2011} (see Figure \ref{subfig-doubling-mutated-a3}) admits an involution and is without oriented cycle. 

\subsection{Action of an involution}

In the remainder of Section \ref{section-quiver-with-involutions}, we assume that $Q$ (without oriented cycle) admits an involution $\tau$.
As indicated in the introduction, the involution $\tau$ acts on many objects attached to the quiver $Q$. 
First, recall that $\tau$ yields an involution on $\Z^{Q_0}$ and remark that, for all $\s,\alpha,\beta\in\Z^{Q_0}$,
\begin{align*}
	\pscal{\alpha}{\beta} &= \pscal{\tp \beta}{\tp \alpha}, \\
	(\tp \s)(\beta) &= \s(\tp \beta), \\
	\tp \pscal{\alpha}{\cdot} &= \pscal{\cdot}{\tp\alpha}.
\end{align*}
Let $\alpha\in\N^{Q_0}$.
If $\alpha = \tp \alpha$, $\tau$ is an anti-linear involution on the complex vector space $\Rep_\alpha(Q)$ and we obtain a real structure on it: this yields lemma \ref{lemm-dense-fixed-points-involution}.\\

\begin{lemm} \label{lemm-dense-fixed-points-involution}
	If $\alpha=\tp\alpha$, the set $\Rep_\alpha(Q)^\tau$ of all $\tau$-symmetric representations is Zariski-dense in $\Rep_\alpha(Q)$.
\end{lemm}

\begin{exem}
	We consider the case where $Q$ is $A_2^\rightarrow$ (see Figure \ref{fig-quiver-an}) equiped with its only involution $\tau$ and $\alpha \in \N^{Q_0}$ such that $\alpha=\tp \alpha$.
	The set $\Rep_\alpha(Q)^\tau$ is the Zariski-dense set of all hermitian matrices.
	On the other hand, the set of orthogonal (resp. symplectic) representations defined in \cite[Section 2]{derksen_generalized_2002} is the Zariski-closed (and not dense) set of all symmetric matrices (resp. of all skew-symmetric matrices).
\end{exem}

For all $n\in\N$ and $A\subset \C^{n}$, we denote by $A^\perp$ the orthogonal of $A$ in $\C^n$ with respect to the canonical Hermitian product $(z,z') \mapsto \sum_{i\in[n]} \overline{z_i } z'_i$.
For all representation $V$ and all subrepresentation $W<V$ we define the $\tau$-orthogonal of $W$ by 
$$W^\perp = (W(\tau x)^\perp)_x <\tp V,$$
which is a subrepresentation of $\tp V$ of dimension 
$$\dim W^\perp = \tp(\dim V-\dim W)^.$$

\begin{defi}
	A $\tau$-\textit{isotropic} subrepresentation of a $\tau$-symmetric representation $V$ is a subrepresentation $W<V$ satisfying $W<W^\perp$. 
\end{defi}

More generally, we have a bijective map from the filtrations of $V$ of length $r$ to the filtrations of $\tp V$ of length $r$: to a filtration $(F_i)_{i}$ of dimension $(\alpha_i)_{i}$ we associate $\tp F = (F_{r+1-i}^\perp)$ of dimension $\tp (\alpha_i) = (\tp \alpha_{r+1-i})$. From this we know that $(\alpha_i)_{i}$ is a generic filtration dimension if and only if $\tp (\alpha_i)_{i}$ is a generic filtration dimension. In particular, $\beta\hookrightarrow \alpha$ if and only if $\tp(\alpha-\beta) \hookrightarrow \tp \alpha$.

For all representation $V$ and $W$, $\tau$ defines an anti-linear isomorphism $\Hom_Q(V,W) \rightarrow \Hom_Q(\tp W,\tp V)$ (which is self-inverse if $W= \tp V$) by $\tp \ph = (\ph(\tau x)^*)_x$. From this and Equation \eqref{eq-ringel=hom-ext} we deduce that
\begin{align} \label{eq-dual-generic-hom}
	\hom_Q(\alpha,\beta) &= \hom_Q(\tp\beta,\tp\alpha),\\
	\ext_Q(\alpha,\beta) &= \ext_Q(\tp\beta,\tp\alpha). \label{eq-dual-generic-ext}
\end{align}
Since $\tau$ yields an anti-linear isomorphism $\SI(Q,\alpha)_\s \simeq\SI(Q,\tau\cdot \alpha)_{-\tau \cdot \s}$, by Equation \eqref{eq-reciprocity-property} we have
\begin{equation} \label{eq-involution-alpha-circ-beta}
	\alpha\circ \beta = (\tp\beta)\circ (\tp\alpha).
\end{equation}

\subsection{The minimum witness}

Assume that $\alpha = \tp \alpha$, that $V\in \Rep_\alpha(Q)$ is a $\tau$-symmetric representation and that $\s = -\tp \s$. Let $W$ be the minimum witness of $(V,\s)$, $\beta:= \dim W$ and $\ga:= \alpha-\beta-\tp\beta$.\\

\begin{lemm} \label{lemm-isotropic-witness}
	The maximum witness of $(V,\s)$ is $W^\perp$. In particular, $W$ is $\tau$-isotropic, $\ga\geqslant 0$, $\ga=\tp \ga$ and, for all $x\in Q_0$ such that $\tau x = x$, $\beta(x)=0$.
	In addition to this we have $\hom_Q(\beta,\tp\beta) = \hom_Q(\beta, \ga)=0.$
\end{lemm}

\begin{proof}
	Let $M$ the maximum witness. We have 
	$$\s(\dim W) = \s(\dim W^\perp) = \s(\dim M) = \s(\dim M^\perp)$$
	hence $W^\perp<M$ and $W<M^\perp$, i.e. $M=W^\perp$ and $W<W^\perp$.
	
	Below, the kernel (resp. image) of a morphism $\ph$ will be denoted by $K$ (resp. $I$) and its dimension $k:= \dim K$ (resp. $i:= \dim I$).
	
	Let $\ph\in\Hom_Q(W,V/W^\perp)$, $\pi:V\twoheadrightarrow V/W^\perp$ the canonical projection and $\tilde I:= \pi^{-1}(I)$. 
	Then $\tilde I$ is a subrepresentation of $V$ of dimension $\dim I + \dim W^\perp$ and we have
	$$\disc(V,\s) \geqslant \s(\dim \tilde I) = \s(i) + \s(\dim W^\perp) = \disc(V,\s)+\s(i).$$
	Then $\s(i)\leqslant 0$. Yet, by the rank theorem, $\beta = i+k$ and $$\disc(V,\s)= \s(\beta) = \s(i)+\s(k)\leqslant \s(k)$$
	hence $\disc(V,\s)\leqslant \s(k)$. Then the kernel $K$ is a witness for $(V,\s)$ and a subrepresentation of the minimum witness $W$, i.e. $K=W$, i.e. $\ph=0$. From this we deduce the first equality $\hom_Q(\beta,\tp \beta) = 0.$
	
	Let $\ph\in\Hom_Q(W,W^\perp/W)$, $\pi:W^\perp\twoheadrightarrow W^\perp/W$ the canonical projection and $\tilde I:= \pi^{-1}(I)$. Then $\tilde I$ is a subrepresentation of $W^\perp<V$ of dimension $\dim I + \dim W$ and we have
	$$\disc(V,\s) \geqslant \s(\dim \tilde I) = \s(i)+\disc(V,\s)$$
	hence $\s(i)\leqslant 0$. As before, by the rank theorem, we have 
	$$\disc(V,\s)= \s(\beta) = \s(i)+\s(k)\leqslant \s(k).$$
	As before, from this we deduce $K=W$, i.e. $\ph=0$ and $\hom_Q(\beta,\ga)=0$.
\end{proof}

Note that Schofield's Equation \eqref{eq-schofield} shows that $\hom_Q(\cdot,\cdot)$ and $\ext_Q(\cdot,\cdot)$ are subadditive with respect to each variable. From this, Lemma \ref{lemm-isotropic-witness} and Equation \eqref{eq-dual-generic-hom}, we have
$$\hom_Q(\ga,\tp\beta) = \hom_Q(\beta,\tp\beta+\ga) = \hom_Q(\ga+\beta,\tp\beta) =0.$$

\subsection{The antisymmetric elements} \label{section-proof-main-result}

Let $\alpha\in \N^{Q_0}$ such that $\alpha = \tp\alpha$. Recall that $\Sigma(Q,\alpha)$ is stable by $\tau$.
We study the subsemigroup $\Sigma(Q,\alpha)^{-\tau}$ of all $\s\in \Sigma(Q,\alpha)$ such that $\s=-\tp \s$. \\

\begin{defi}
	Let $\II^{00}(Q,\alpha,\tau)$ be the set of all $(\beta,\ga) \in (\N^{Q_0})^2$ such that $\alpha = \beta + \ga + \tp \beta$ and $\beta \circ \ga = \beta \circ \tp \beta = 1$.
\end{defi}

For all $(\beta,\ga), (\beta',\ga') \in \II^{00}(Q,\alpha,\tau)$, $\beta = \beta' \Rightarrow \ga = \ga'$; considering pairs of dimensions is simply a way to avoid redefining $\ga$ each time.
Note that $\II^{00}(Q,\alpha,\tau)$ depends on $\tau_0$ but not on $\tau_1$.
By Equation \eqref{eq-involution-alpha-circ-beta} and Proposition \ref{prop-ressayre-multiplicative-coro2}, for all $(\beta,\ga) \in \II^{00}(Q,\alpha,\tau)$,
\begin{equation*}
	\ga \circ \tp\beta = \beta \circ (\ga + \tp\beta) = (\ga +\beta)\circ \tp\beta = 1
\end{equation*}
and, in particular: $(\beta,\ga,\tp\beta)$ is a generic filtration dimension, $\beta \hookrightarrow\alpha$ and $\pscal{\beta}{\cdot}-\pscal{\cdot}{\tp\beta} \in \Sigma(Q,\ga)^{-\tau}$.
This set of dimensions is important in the following lemma.\\

\begin{lemm} \label{lemm-disc-calculation-involution-best}
	For all $\s\in \Z^{Q_0}$ such that $\s= -\tp \s$ we have
	$$\disc(\alpha,\s) = \max_{ (\beta,\ga) \in \II^{00}(Q,\alpha,\tau) } \s(\beta) = \max_{ (\beta,\ga) \in \II^{00}(Q,\alpha,\tau), \ga \in \overline{\Sigma}(Q,\s) } \s(\beta).$$
\end{lemm}

\begin{proof}
	Derksen-Weyman's Theorem \ref{theo-dw-sigma-is-a-cone} yields
	$$\disc(\alpha,\s) \geqslant \max_{ (\beta,\ga) \in \II^{00}(Q,\alpha,\tau), \ga \in \overline{\Sigma}(Q,\s) } \s(\beta)$$
	We will now prove the reversed inequality.
	By Lemma \ref{lemm-filtration-zariski-closed}, there exists a non-empty open subset $\UU$ of $\Rep_\alpha(Q)$ such that, for all $V\in \UU$ and all $\beta \in \N^{Q_0}$, 
	\begin{enumerate}
		\item if $V$ admits a subrepresentation of dimension $\beta$ then $\beta \hookrightarrow \alpha$;
		\item if $V$ admits a filtration of dimension $(\beta,\alpha-\beta-\tp\beta,\tp\beta)$ then it is a generic filtration dimension.
	\end{enumerate}
	We denote by $\II(Q,\alpha,\tau)$ the set of all $(\beta,\ga)\in (\N^{Q_0})^2$ such that
	$\alpha = \beta + \ga+\tp\beta$ and $\pscal{\beta}{\ga}= \pscal{\beta}{\tp\beta} = 0$.
	By Derksen-Weyman-Schofield's Theorem \ref{theo-dw-schofield-number-of-subrepresentations}, for all $(\beta,\ga)\in \II(Q,\alpha,\tau)$, there exists a non-empty open subset $\VV_{\beta}\subset \Rep_\alpha(Q)$ such that any element of $\VV_{\beta}$ admits exactly $\beta\circ \ga + \tp\beta$ subrepresentations of dimension $\beta$. 
	By Corollary \ref{coro-general-sub-sub-representations}, for all $(\beta,\ga)\in \II(Q,\alpha,\tau)$, there exists a non-empty open subset $\WW_{\beta}\subset \Rep_\alpha(Q)$ such that, for all $V\in \WW_\beta$ and all $(\beta+\ga)$-dimensional subrepresentation $W<V$, $W$ admits exactly $\beta\circ \ga$ subrepresentations of dimension $\beta$. Let 
	\begin{align*}
		\VV &= \bigcap_{(\beta,\ga) \in \II(Q,\alpha,\tau)} \VV_\beta, &
		\WW = \bigcap_{(\beta,\ga) \in \II(Q,\alpha,\tau)} \WW_\beta.
	\end{align*} 
	By Lemma \ref{lemm-dense-fixed-points-involution}, there exists a $\tau$-symmetric representation $V\in \UU\cap \VV\cap \WW$. 
	Let $\s\in \Z^{Q_0}$ such that $\s = -\tp \s$. Let $W$ be the minimum witness of $(V,\s)$, $\beta = \dim W$ its dimension and $\ga = \alpha-\beta-\tp\beta$. 
	Since $V\in\UU$, $$\s(\beta) = \disc(\alpha,\s).$$
	By Lemma \ref{lemm-isotropic-witness}, the subrepresentation $W$ is $\tau$-isotropic and $V$ admits a filtration of dimension $(\beta, \ga, \tp\beta)$. Since $V\in\UU$, $(\beta,\ga,\tp\beta)$ is a generic filtration dimension.
	Then Lemma \ref{lemm-sum-ringel-generic-dim-vector} yields 
	$$2\pscal{\beta}{\ga} + \pscal{\beta}{\tp\beta}\geqslant 0.$$ 
	Yet, by Equation \eqref{eq-ringel=hom-ext} and Lemma \ref{lemm-isotropic-witness}, $$\pscal{\beta}{\ga},\pscal{\beta}{\tp\beta} \leqslant 0.$$ 
	From this we know that $(\beta,\ga)\in \II(Q,\alpha,\tau)$.
	
	Since $V\in\VV$ then, by Lemma \ref{lemm-unique-witness}, $\beta\circ \ga+\tp\beta = 1$.
	Since $V\in\WW$ then, by Lemma \ref{lemm-unique-witness} and Lemma \ref{lemm-isotropic-witness}, $\beta\circ \ga=1$, i.e. $\ga\circ \tp\beta =1$.
	Then, by Corollary \ref{prop-ressayre-multiplicative-coro2}, $\beta\circ \tp\beta =1$ and $(\beta,\ga)\in \II^{00}(Q,\alpha,\tau)$.
	
	It remains to be proven that $\s\in\Sigma(Q,\ga)$. We consider the quotient $W^\perp/W$ and $\pi:W^\perp \twoheadrightarrow W^\perp/W$ the canonical projection. Let $X$ be a subrepresentation of $W^\perp/W$ and $\tilde X:= \pi^{-1}(X) < W^\perp <V$ of dimension $\dim X + \beta$. 
	Since $\tilde X$ is a subrepresentation of $V$, $\s(\dim \tilde X) \leqslant \s(\beta)$, i.e. $\s(X)\leqslant 0$. By King's criterion \ref{theo-king-criterion}, $W^\perp/W$ is $\s$-semi-stable. By Lemma \ref{lemm-sigma-made-of-semi-stable}, this shows that $\s\in\Sigma(Q,\ga)$. 
\end{proof}

Lemma \ref{lemm-disc-calculation-involution-best} is comparable to Lemma \ref{lemm-disc-calculation-circ1}.
Associated with Derksen-Weyman's \ref{theo-dw-sigma-is-a-cone}, it yields Theorem \ref{theo-symmetric-sigma-inequalities-circ1} below, which is the main result of this paper and yields Theorem \ref{theo-coro-main}.\\

\begin{theo} \label{theo-symmetric-sigma-inequalities-circ1}
	Let $\alpha\in \N^{Q_0}$ and $\s\in\Z^{Q_0}$ such that $\alpha = \tp \alpha$ and $\s= -\tp\s$.
	We have $\s\in \Sigma(Q,\alpha)^{-\tau}$ if and only if $\s(\beta)\leqslant 0$ for all $\beta\in\N^{Q_0}$ such that
	\begin{enumerate}
		\item $\ga:= \alpha-\beta-\tp\beta\in \N^{Q_0}$;
		\item $\beta \circ \ga = \beta \circ \tp \beta = 1$;
		\item $\sigma\in\Sigma(Q,\ga)^{-\tau}$.
	\end{enumerate}
\end{theo}

\section{Examples} \label{section-examples}

\subsection{Notations} 

Let $Q$ be a quiver, $\tau$ an involution (see Definition \ref{defi-involution}) on $Q$ and $\alpha$. We denote by $\II^0(Q,\alpha,\tau)$ the set of all $(\beta,\ga) \in(\N^{Q_0})^2$ such that $\alpha = \beta + \ga+ \tp \beta$, $\beta \circ \ga \neq 0$ and $\beta \circ \tp\beta \neq 0.$
There is an algorithm to compute the following integers:
\begin{align*}
	n_1(Q,\alpha) &= \Card \ens{ \beta \leqslant \alpha | \beta \hookrightarrow \alpha },\\
	n_2(Q,\alpha) &= \Card \ens{ \beta \leqslant \alpha | \beta \circ(\alpha - \beta) \neq 0 },\\
	n_3(Q,\alpha,\tau) &= \Card \II^0(Q,\alpha,\tau).
\end{align*}
Remark that $n_1(Q,\alpha)$ (resp. $n_2(Q,\alpha)$) is the number of inequalities given by Derksen-Weyman's Theorem \ref{theo-dw-sigma-is-a-cone} (resp. Theorem \ref{theo-inductive-sigma}) to characterise $\Sigma(Q,\alpha)$ and $n_3(Q,\alpha,\tau)$ is the number of inequalities given by Theorem \ref{theo-coro-main} to characterise $\Sigma(Q,\alpha)^{-\tau}$.

\subsection{A quiver of type $\hat D_5$} Figure \ref{fig-number-of-ineq-D5} gives some values of the maps $n_1$, $n_2$ and $n_3$ for the quiver and the involution defined in example \ref{exemquiver1}. A dimension $\alpha \in \N^{Q_0}$ is denoted by $(\alpha(x_1))_{i\in [6]}$: then 
$$\ens{\alpha\in \N^{Q_0} | \alpha= \tp\alpha} = \ens{ (a,b,c,c,b,a) ; a,b,c\in\N }.$$

\begin{figure}[htbp]
	\centering
	\caption{Numbers of inequalities for Example \ref{exemquiver1}.}
	\label{fig-number-of-ineq-D5} 
	\begin{tabular}{|c|c|c|c|} 
		\hline
		$\alpha$           & $n_1(Q,\alpha)$ & $n_2(Q,\alpha)$ & $n_3(Q,\alpha,\tau)$ \\\hline
		$(1,1,1,1,1,1)$ & $9$        & $9$        & $5$           \\\hline
		$(2,2,2,2,2,2)$ & $43$       & $9$        & $5$           \\\hline
		$(3,3,3,3,3,3)$ & $147$      & $9$        & $5$           \\\hline
		$(4,4,4,4,4,4)$ & $406$      & $9$        & $5$           \\\hline
		$(0,1,2,2,1,0)$ & $12$       & $9$        & $4$           \\\hline
		$(0,2,1,1,2,0)$ & $16$       & $12$       & $6$            \\\hline
		$(2,1,0,0,1,2)$ & $36$       & $16$       & $9$           \\\hline
		$(1,2,3,3,2,1)$ & $59$       & $25$       & $7$           \\\hline
		$(2,3,2,2,3,2)$ & $112$      & $20$       & $9$           \\\hline
		$(2,3,4,4,3,2)$ & $244$      & $57$       & $10$\\
		\hline              
	\end{tabular}
\end{figure}

\subsection{The Sun quiver} \label{subsection-sun-quiver}
\begin{figure}[htbp]
	\centering
	\caption{Involutions on the $(6,1)$-Sun quiver.}
	\label{fig-involutions-sun-quiver}
	\begin{subfigure}{0.23\textwidth}
		\caption{Involution $\tau$} \label{subfig-tau-involution}
		\[\begin{tikzcd}
			5 & 4 \\
			6 & 3 \\
			1 & 2
			\arrow[color={red}, dotted, no head, from=1-1, to=3-2]
			\arrow[from=1-2, to=1-1]
			\arrow[from=1-2, to=2-2]
			\arrow[from=2-1, to=1-1]
			\arrow[from=2-1, to=3-1]
			\arrow[shift right, color={red}, dotted, no head, from=2-1, to=3-1]
			\arrow[shift right, color={red}, dotted, no head, from=2-2, to=1-2]
			\arrow[from=3-2, to=2-2]
			\arrow[from=3-2, to=3-1]
		\end{tikzcd}\]
	\end{subfigure}
	\hspace{1cm}
	\begin{subfigure}{0.23\textwidth}
		\caption{Involution $\rho$} \label{subfig-rho-involution}
		\[\begin{tikzcd}
			5 & 4 \\
			6 & 3 \\
			1 & 2
			\arrow[from=1-2, to=1-1]
			\arrow[from=1-2, to=2-2]
			\arrow[from=2-1, to=1-1]
			\arrow[from=2-1, to=3-1]
			\arrow[color={red}, dotted, no head, from=2-2, to=2-1]
			\arrow[color={red}, dotted, no head, from=3-1, to=1-2]
			\arrow[color={red}, dotted, no head, from=3-2, to=1-1]
			\arrow[from=3-2, to=2-2]
			\arrow[from=3-2, to=3-1]
		\end{tikzcd}\]
	\end{subfigure}
\end{figure}
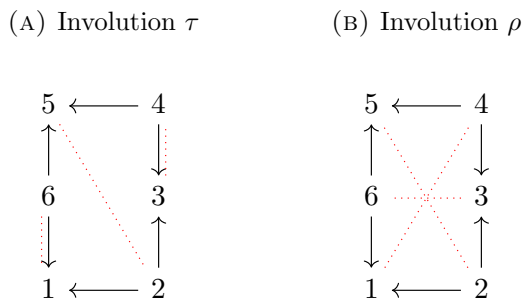

Here we define the Sun quiver introduced in \cite[Section 3.2]{collins_generalized_2020}. 
Let $k\geqslant 2$ and $n\geqslant 1$. Let $Q_0= \Z/2k\Z \x [n]$ be a set of vertices and $Q_1= \ens{a_{i,j}; (i,j)\in Q_0 }$ be a set of arrows. Since $2k$ is even, we can talk about the parity of an element in $\Z/2k\Z$. We define $t,h:Q_1\rightarrow Q_0$ by: for all $i\in \Z/2k$, if $i$ is odd then
\begin{align*}
	ta_{i,n} &= (i+1,n),&
	ha_{i,n} &= (i,n),\\
	\forall j \in [n-1], ta_{i,j} &= (i,j+1), &
	ha_{i,j} &= (i,j),		
\end{align*}
and if $i$ is even then
\begin{align*}
	ta_{i,n} &= (i,n),&
	ha_{i,n} &= (i+1,n),\\
	\forall j \in [n-1], ta_{i,j} &= (i,j), &
	ha_{i,j} &= (i,j+1).		
\end{align*}
Then $(Q_0,Q_1,h,t)$ is called the $(2k,n)$-\textit{Sun quiver}.
Let $\ta_0:Q_0\rightarrow Q_0$ and $\ta_1:Q_1\rightarrow Q_1$ the only self-inverse maps defined by: for all $i\in \range{0,k}$ and $j\in[n]$, $\ta_{0}(i,j) = (1-i,j)$, $\ta_{1} a_{i,n} = a_{-i,n}$ and, if $j\leqslant n-1$, $\ta_{1} a_{i,n} = a_{1-i,j}$.
The couple $(\ta_{0}, \ta_{1})$ defines the unique involution $\ta$ such that $\ta (0,n) = (1,n)$. See Figure \ref{subfig-tau-involution}. 
The two maps
$$ \appli{\rho_0}{ (i,j) }{Q_0}{ (i+k,j) }{Q_0}, \appli{\rho_1}{ a_{i,j} }{Q_1}{ a_{i+k,j} }{Q_1}$$
define an involution $\rho$ on $Q$ if and only if $k$ is odd: see Figure \ref{subfig-rho-involution}.

We consider the case where $(k,n) = (3,1)$, which is studied in \cite[Section 9]{baldoni_horn_2023}.
Figure \ref{fig-number-of-ineq-sun} gives some values of the integers $n_1$, $n_2$ and $n_3$ for the $(6,1)$-Sun quiver and the two involutions previously defined. A dimension $\alpha \in \N^{Q_0}$ is denoted by $(\alpha(i,1))_{i\in [6]}$.
With this notation we have
\begin{align*}
	\ens{\alpha\in \N^{Q_0} | \alpha= \tp\alpha} &= \ens{ (a,a,b,c,c,b) ; a,b,c\in\N }, \\
	\ens{\alpha\in \N^{Q_0} | \alpha= \rho\cdot \alpha} &= \ens{ (a,b,c,a,b,c) ; a,b,c\in\N }, \\
	\ens{\alpha\in \N^{Q_0} | \alpha= \tp\alpha = \rho\cdot \alpha} &= \ens{ (a,a,a,a,a,a) ; a\in\N }.
\end{align*}

\begin{figure}[htbp]
	\centering
	\caption{Numbers of inequalities for the $(6,1)$-Sun quiver.}
	\label{fig-number-of-ineq-sun}
	\begin{tabular}{|c|c|c|c|c|} 
		\hline
		$\alpha$           & $n_1(Q,\alpha)$ & $n_2(Q,\alpha)$ & $n_3(Q,\alpha,\tau)$ & $n_3(Q,\alpha,\rho)$  \\\hline
		$ (1,1,2,3,3,2)$& $159$ & $48$ & $12$ & \\\hline
		$ (1,1,3,2,2,3)$& $155$ & $32$ & $11$ & \\\hline
		$ (2,2,2,4,4,2)$& $396$ & $30$ & $11$ & \\\hline
		$ (2,2,3,4,4,3)$& $717$ & $54$ & $14$ & \\\hline
		$ (1,2,3,1,2,3)$& $190$ & $56$ &      & $25$ \\\hline
		$ (2,1,3,2,1,3)$& $190$ & $56$ &      & $25$ \\\hline
		$ (2,4,2,2,4,2)$& $464$ & $32$ &      & $17$ \\\hline
		$ (3,2,4,3,2,4)$& $924$ & $112$ &      & $32$ \\\hline
		$ (2,2,2,2,2,2)$& $129$ & $19$ &  $10$    & $12$ \\\hline
		$ (3,3,3,3,3,3)$& $571$ & $20$ &  $10$    & $12$ \\\hline          
	\end{tabular}
\end{figure}

\bibliography{bibtemplate}

@article{baldoni_horn_2023,
	title = {Horn conditions for quiver subrepresentations and the moment map},
	volume = {19},
	issn = {15588599, 15588602},
	doi = {10.4310/PAMQ.2023.v19.n4.a1},
	language = {en},
	number = {4},
	journal = {Pure and Applied Mathematics Quarterly},
	author = {Baldoni, Velleda and Vergne, Michèle and Walter, Michael},
	year = {2023},
	pages = {1687--1731},
}

@article{berline_horn_2018,
	title = {The {Horn} inequalities from a geometric point of view},
	volume = {63},
	issn = {0013-8584, 2309-4672},
	doi = {10.4171/lem/63-3/4-7},
	language = {en},
	number = {3},
	journal = {L’Enseignement Mathématique},
	author = {Berline, Nicole and Vergne, Michèle and Walter, Michael},
	month = sep,
	year = {2018},
	pages = {403--470},
}

@article{derksen_semi-invariants_2000,
	title = {Semi-invariants of quivers and saturation for {Littlewood}-{Richardson} coefficients},
	volume = {13},
	issn = {0894-0347, 1088-6834},
	doi = {10.1090/S0894-0347-00-00331-3},
	language = {en},
	number = {3},
	journal = {Journal of the American Mathematical Society},
	author = {Derksen, Harm and Weyman, Jerzy},
	month = mar,
	year = {2000},
	pages = {467--479},
}

@article{derksen_generalized_2002,
	title = {Generalized quivers associated to reductive groups},
	volume = {94},
	issn = {0010-1354, 1730-6302},
	doi = {10.4064/cm94-2-1},
	language = {en},
	number = {2},
	journal = {Colloquium Mathematicum},
	author = {Derksen, Harm and Weyman, Jerzy},
	year = {2002},
	pages = {151--173},
}

@article{derksen_combinatorics_2011,
	title = {The combinatorics of quiver representations},
	volume = {61},
	issn = {0373-0956, 1777-5310},
	doi = {10.5802/aif.2636},
	language = {en},
	number = {3},
	journal = {Annales de l’institut Fourier},
	author = {Derksen, Harm and Weyman, Jerzy},
	year = {2011},
	pages = {1061--1131},
}

@article{horn_eigenvalues_1962,
	title = {Eigenvalues of sums of {Hermitian} matrices},
	volume = {12},
	issn = {0030-8730, 0030-8730},
	doi = {10.2140/pjm.1962.12.225},
	language = {en},
	number = {1},
	journal = {Pacific Journal of Mathematics},
	author = {Horn, Alfred},
	month = mar,
	year = {1962},
	pages = {225--241},
}

@article{king_moduli_1994,
	title = {Moduli of representations of finite dimensional algebras},
	volume = {45},
	issn = {0033-5606, 1464-3847},
	doi = {10.1093/qmath/45.4.515},
	language = {en},
	number = {4},
	journal = {The Quarterly Journal of Mathematics},
	author = {King, A. D.},
	year = {1994},
	pages = {515--530},
}

@article{schofield_general_1992,
	title = {General {Representations} of {Quivers}},
	volume = {s3-65},
	issn = {00246115},
	doi = {10.1112/plms/s3-65.1.46},
	language = {en},
	number = {1},
	journal = {Proceedings of the London Mathematical Society},
	author = {Schofield, Aidan},
	month = jul,
	year = {1992},
	pages = {46--64},
}

@book{derksen_introduction_2017,
	series = {Graduate {Studies} in {Mathematics}},
	title = {An {Introduction} to {Quiver} {Representations}},
	volume = {184},
	isbn = {978-1-4704-2556-2 978-1-4704-4260-6},
	language = {en},
	publisher = {American Mathematical Society},
	author = {Derksen, Harm and Weyman, Jerzy},
	month = nov,
	year = {2017},
	doi = {10.1090/gsm/184},
}

@article{fulton_expose_1998,
	series = {Astérisque},
	title = {Exposé {Bourbaki} 845 : {Eigenvalues} of sums of {Hermitian} matrices [after {A}. {Klyachko}]},
	volume = {1997/98},
	issn = {03031179, 24925926},
	shorttitle = {Exposé {Bourbaki} 845},
	doi = {10.24033/ast.426},
	language = {en},
	number = {252},
	journal = {Astérisque},
	author = {Fulton, William},
	year = {1998},
	pages = {255--269},
}

@misc{huszar_non-commutative_2021,
	title = {Non-commutative {Rank} and {Semi}-stability of {Quiver} {Representations}},
	language = {en},
	publisher = {arXiv},
	author = {Huszar, Alana},
	month = oct,
	year = {2021},
	note = {arXiv:2111.00039 [math]},
}

@article{ressayre_multiplicative_2011,
	title = {Multiplicative formulas in {Schubert} calculus and quiver representation},
	volume = {22},
	issn = {0019-3577},
	doi = {https://doi.org/10.1016/j.indag.2011.08.004},
	language = {en},
	number = {1},
	journal = {Indagationes Mathematicae},
	author = {Ressayre, Nicolas},
	year = {2011},
	pages = {87--102},
}

@article{schofield_semi-invariants_1991,
	title = {Semi-{Invariants} of {Quivers}},
	volume = {s2-43},
	issn = {0024-6107},
	doi = {10.1112/jlms/s2-43.3.385},
	language = {en},
	number = {3},
	journal = {Journal of the London Mathematical Society},
	author = {Schofield, Aidan},
	month = jun,
	year = {1991},
	pages = {385--395},
}

@article{collins_generalized_2020,
	title = {Generalized {Littlewood}–{Richardson} coefficients for branching rules of {GL}(n) and extremal weight crystals},
	volume = {3},
	issn = {2589-5486},
	doi = {10.5802/alco.143},
	language = {en},
	number = {6},
	journal = {Algebraic Combinatorics},
	author = {Collins, Brett},
	month = dec,
	year = {2020},
	pages = {1365--1400},
}

@article{ressayre_git-cones_2012,
	title = {{GIT}-cones and quivers},
	volume = {270},
	issn = {0025-5874, 1432-1823},
	doi = {10.1007/s00209-010-0796-0},
	language = {en},
	number = {1-2},
	journal = {Mathematische Zeitschrift},
	author = {Ressayre, Nicolas},
	month = feb,
	year = {2012},
	pages = {263--275},
}

@article{chindris_quivers_2009,
	title = {Quivers, long exact sequences and {Horn} type inequalities {II}},
	volume = {51},
	issn = {0017-0895, 1469-509X},
	doi = {10.1017/S0017089508004631},
	language = {en},
	number = {2},
	journal = {Glasgow Mathematical Journal},
	author = {Chindris, Calin},
	month = may,
	year = {2009},
	pages = {201--217},
}

@article{bialynicki-birula_theorems_1973,
	title = {Some {Theorems} on {Actions} of {Algebraic} {Groups}},
	volume = {98},
	issn = {0003486X},
	doi = {10.2307/1970915},
	language = {en},
	number = {3},
	journal = {The Annals of Mathematics},
	author = {Białynicki-Birula, A.},
	month = nov,
	year = {1973},
	pages = {480--497},
}

@article{chindris_simultaneous_2021,
	title = {Simultaneous robust subspace recovery and semi-stability of quiver representations},
	volume = {577},
	issn = {00218693},
	doi = {10.1016/j.jalgebra.2021.03.005},
	language = {en},
	journal = {Journal of Algebra},
	author = {Chindris, Calin and Kline, Daniel},
	month = jul,
	year = {2021},
	pages = {210--236},
}

@misc{fomin_total_2010,
	title = {Total positivity and cluster algebras},
	doi = {10.48550/arXiv.1005.1086},
	language = {en},
	publisher = {arXiv},
	author = {Fomin, Sergey},
	month = may,
	year = {2010},
	note = {arXiv:1005.1086 [math]},
}

@article{derksen_number_2007,
	title = {On the number of subrepresentations of a general quiver representation},
	volume = {76},
	issn = {00246107},
	doi = {10.1112/jlms/jdm043},
	language = {en},
	number = {1},
	journal = {Journal of the London Mathematical Society},
	author = {Derksen, Harm and Schofield, Aidan and Weyman, Jerzy},
	month = aug,
	year = {2007},
	pages = {135--147},
}

@article{belkale_local_2001,
	title = {Local systems on {P1}-{S} for {S} a finite set},
	volume = {129},
	issn = {0010437X},
	doi = {10.1023/A:1013195625868},
	language = {en},
	number = {1},
	journal = {Compositio Mathematica},
	author = {Belkale, Prakash},
	year = {2001},
	pages = {67--86},
}

@article{belkale_geometric_2006,
	title = {Geometric proofs of {Horn} and saturation conjectures},
	volume = {15},
	issn = {1056-3911, 1534-7486},
	doi = {10.1090/S1056-3911-05-00420-0},
	language = {en},
	number = {1},
	journal = {Journal of Algebraic Geometry},
	author = {Belkale, Prakash},
	year = {2006},
	pages = {133--173},
}

@article{knutson_honeycomb_1999,
	title = {The honeycomb model of {$GL_n(\mathbb C)$} tensor products {I}: Proof of the saturation conjecture},
	volume = {12},
	issn = {0894-0347, 1088-6834},
	shorttitle = {The honeycomb model of \${GL}\_n({\textbackslash}mathbb {C})\$ tensor products {I}},
	doi = {10.1090/S0894-0347-99-00299-4},
	language = {en},
	number = {4},
	journal = {Journal of the American Mathematical Society},
	author = {Knutson, Allen and Tao, Terence},
	year = {1999},
	pages = {1055--1090},
}
\bibliographystyle{smfalpha}
\end{document}